\documentclass[12pt]{article}
\usepackage{latexsym,amsmath}
\usepackage{amssymb}

\newtheorem{thm}{Theorem}[section]

\newtheorem{lemma}[thm]{Lemma}

\newtheorem{cor}[thm]{Corollary}
\newtheorem{remark}[thm]{Remark}

\newtheorem{claim}[thm]{Claim}

\newcommand{\beq}[1]{\begin{equation}\label{#1}}
\newcommand{\enq}[0]{\end{equation}}

\newcommand{\qed}[0]{{\hspace*{\fill}\mbox{$\Box$}}}

\newcommand{\cA}[0]{{\cal A}}
\newcommand{\cE}[0]{{\cal E}}
\newcommand{\cI}[0]{{\cal I}}
\newcommand{\cJ}[0]{{\cal J}}
\newcommand{\cM}[0]{{\cal M}}
\newcommand{\cO}[0]{{\cal O}}

\newcommand{\cU}[0]{{\cal U}}
\newcommand{\cV}[0]{{\cal V}}
\newcommand{\cW}[0]{{\cal W}}

\newcommand{\Z}{{\mathbb Z}}

\newcommand{\ga}[0]{\alpha}
\newcommand{\gb}[0]{\beta}
\newcommand{\gre}[0]{\epsilon}
\newcommand{\gd}[0]{\delta}
\newcommand{\gD}[0]{\Delta}
\newcommand{\grg}[0]{\gamma}
\newcommand{\gG}[0]{\Gamma}
\newcommand{\gk}[0]{\kappa}
\newcommand{\gl}[0]{\lambda}
\newcommand{\go}[0]{\omega}
\newcommand{\gO}[0]{\Omega}
\newcommand{\gS}[0]{\Sigma}

\begin{document}

\renewcommand{\thefootnote}{\fnsymbol{footnote}}

\title{Slow mixing of Glauber Dynamics for the hard-core model on regular
bipartite graphs}

\author{David Galvin\thanks{Institute for Advanced Study, Einstein Drive, Princeton, NJ 08540; galvin@ias.edu.
Research supported in part by NSF grant DMS-0111298.}\\
Prasad Tetali\thanks{School of Mathematics \& College of
Computing, Georgia Institute of Technology, Atlanta, GA
30332-0160. Research supported in part by NSF grant DMS-0100289.}}

\date{Appeared 2006}

\maketitle

\begin{abstract}
Let $\gS=(V,E)$ be a finite, $d$-regular bipartite graph. For any
$\gl>0$ let $\pi_\gl$ be the probability measure on the
independent sets of $\gS$ in which the set $I$ is chosen with
probability proportional to $\gl^{|I|}$ ($\pi_\gl$ is the {\em
hard-core measure with activity $\gl$ on $\gS$}). We study the
Glauber dynamics, or single-site update Markov chain, whose
stationary distribution is $\pi_\gl$. We show that when $\gl$ is
large enough (as a function of $d$ and the expansion of subsets of
single-parity of $V$) then the convergence to stationarity is
exponentially slow in $|V(\gS)|$. In particular, if $\gS$ is the
$d$-dimensional hypercube $\{0,1\}^d$ we show that for values of
$\gl$ tending to $0$ as $d$ grows, the convergence to stationarity
is exponentially slow in the volume of the cube. The proof
combines a conductance argument with combinatorial enumeration
methods.
\end{abstract}

\section{Introduction and statement of the result}
\label{sec-intro}

Let $\gS=(V,E)$ be a simple, loopless, finite graph on vertex set
$V$ and edge set $E$. (For graph theory basics, see e.g.
\cite{Bollobas}, \cite{Diestel}.) Write $\cI(\gS)$ for the set of
independent sets (sets of vertices spanning no edges) in $V$. For
$\gl > 0$ we define the {\em hard-core measure with activity
$\gl$} on $\cI(\gS)$ by
\begin{equation} \label{eq-hcmeasure}
\pi_{\gl}(\{I\}) = \frac{\gl^{|I|}}{Z_{\gl}(\gS)} ~~~\mbox{for
$I\in \cI$}
\end{equation}
where $Z_{\gl}(\gS) = \sum_{I\in \cI}\gl^{|I|}$ is the appropriate
normalizing constant. We will often write $w_{\gl}(I)$ for
$\gl^{|I|}$ and, for $\cJ \subseteq \cI$, $w_{\gl}(\cJ)$ for
$\sum_{J \in \cJ} w_{\gl}(J)$.

The hard-core measure originally arose in statistical physics (see
e.g. \cite{Dobrushin,BS}) where it serves as a simple mathematical
model of a gas with particles of non-negligible size. The vertices
of $\gS$ we think of as sites that may or may not be occupied by
particles; the rule of occupation is that adjacent sites may not
be simultaneously occupied. The activity parameter $\gl$ measures
the likelihood of a site being occupied.

The measure also has a natural interpretation in the context of
communications networks (see e.g. \cite{Kelly}). Here the vertices
of $\gS$ are thought of as locations from which ``calls'' can be
made; when a call is made, the call location is connected to all
its neighbours, and throughout its duration, no call may be placed
from any of the neighbours. Thus at any given time, the collection
of locations from which calls are being made is exactly an
independent set in the graph. If calls are attempted independently
at each vertex as a Poisson process of rate $\lambda$ and have
independent exponential mean $1$ lengths, it can be shown that the
long-run stationary distribution of this process is the hard-core
measure on $\gS$.

Our particular focus in this paper is the mixing time of the
Glauber dynamics, or single-site update Markov chain, for this
model. The measure $\pi_{\gl}$ can be realized as the stationary
distribution of a certain Markov chain. Specifically, consider the
chain $\cM_{\gl}=\cM_{\gl}(\gS)$ on state space $\cI(\gS)$ with
transition probabilities $P_{\gl}(I,J), I,J \in \cI(\gS),$ given
by
$$
P_{\gl}(I,J) = \left\{
            \begin{array}{ll}
               0 & \mbox{ if $|I \bigtriangleup J| > 1$} \\
               \frac{1}{|V|}\frac{\gl}{1+\gl}& \mbox{ if $|I \bigtriangleup
J| = 1, ~I \subseteq
               J$} \\
               \frac{1}{|V|}\frac{1}{1+\gl}& \mbox{ if $|I \bigtriangleup J|
= 1, ~J \subseteq
               I$} \\
               1 - \sum_{I \neq J' \in \cI(\gS)} P_{\gl}(I,J') & \mbox{ if
$I=J$}
            \end{array}
         \right.
$$
Underpinning the definition of $\cM_{\gl}$ is the following
dynamical process, known as the {\em Glauber dynamics} on
$\cI(\gS)$. From an independent set $I$, the process follows three
steps. The first step is to choose a vertex $v$ uniformly from
$V$. The second step is to ``add'' $v$ to $I$ with probability
$\gl/(1+\gl)$, and ``remove'' it with probability $1/(1+\gl)$;
that is, to set
$$
I' = \left\{
         \begin{array}{ll}
             I \cup \{v\} & \mbox{ with probability $\frac{\gl}{1+\gl}$} \\
             I \setminus \{v\} & \mbox{ with probability
             $\frac{1}{1+\gl}$}.
         \end{array}
     \right.
$$
The third step is to move to $I'$ if it is a valid independent
set, and stay at $I$ otherwise.

It is readily checked that $\cM_{\gl}$ is an ergodic, aperiodic,
time reversible Markov chain with (unique) stationary distribution
$\pi_{\gl}$.  A natural question to ask about $\cM_{\gl}$ is how
quickly it converges to its stationary distribution. To make this
question precise, we need a few definitions.

Let $\cM$ be an ergodic Markov chain on state space $\gO$, with
transition probabilities $P:\gO^2\rightarrow[0,1]$. For a state
$\go_0 \in \gO$, denote by $P^t(\go_0, \cdot)$ the distribution of
the state at time $t$, given that the initial state is $\go_0$,
and denote by $\pi$ the stationary distribution. Define the {\em
mixing time} of $\cM$ by
$$
\tau_{\cM}=\max_{\go_0 \in \gO} \min
\left\{t_0~:~\frac{1}{2}\sum_{\go \in
\gO}|P^t(\go_0,\go)-\pi(\go)| \leq \frac{1}{e} ~~~ \forall
t>t_0\right\}.
$$
The mixing time of $\cM$ captures the speed at which the chain
converges to its stationary distribution: for every $\gre
>0$, in order to get a sample from $\gO$ which is within $\gre$ of
$\pi$ (in variation distance), it is necessary and sufficient to
run the chain from some arbitrarily chosen distribution for some
multiple (depending on $\gre$) of the mixing time.

\medskip

Much work has been done on the question of bounding
$\tau_{\cM_\gl}$. The strongest general result available to date
is due to Vigoda \cite{Vigoda} who showed that if $\gS$ is any
$n$-vertex graph with maximum degree $\gD$, then
$\tau_{\cM_\gl}(\gS) = O(n\log n)$ whenever $\gl < 2/(\gD-2)$. In
the other direction, Dyer, Frieze and Jerrum
\cite{DyerFriezeJerrum} considered the case $\gl=1$ and showed
that for each $\gD\geq 6$ a random (uniform) $\gD$-regular,
$n$-vertex bipartite $\gS$ almost surely (with probability tending
to $1$ as $n$ tends to infinity) satisfies $\tau_{\cM_1}(\gS) \geq
2^{\gamma n}$ for some absolute constant $\gamma >0$.

Here, we continue in the spirit of \cite{DyerFriezeJerrum} and
construct {\em explicit} families of graphs for which Glauber
dynamics mixes slowly. Specifically, we establish a certain
expansion condition in a regular bipartite graph $\gS$ that forces
$\tau_{\cM_\gl(\gS)}$ to be (almost) exponential in $|V|$ provided
$\gl$ is suitably large (as a function of the expansion). The
$d$-dimensional hypercube $\{0,1\}^d$ satisfies this condition for
$\gl \geq \omega(d^{-1/4}\log^{3/2} d)$.

\medskip

Our work is partly motivated by \cite{Million} where a study was
made of Glauber dynamics for the hard-core measure on the even
discrete torus $T_{L,d}$. This is the graph on
$\{0,\ldots,L-1\}^d$ (with $L$ even) in which two strings are
adjacent if they differ on only one coordinate, and differ by $1
~(\mbox{mod } L)$ on that coordinate. It was shown in
\cite{Million} that for $\gl$ growing exponentially in $d$ (with a
suitably large base), $\tau_{\cM_\gl(T_{L,d})}$ is exponential in
$cL^{d-1}/\log^2 L$ for some $c$ that depends on $d$ but not on
$L$.

In light of a recent result of Galvin and Kahn \cite{GK}, we found
it tempting to believe that slow mixing on $T_{L,d}$ should hold
for much smaller values of $\gl$; even for some values of $\gl$
tending to $0$ as $d$ grows. The main result of \cite{GK} is that
the hard-core model on $\Z^d$ exhibits multiple Gibbs phases for
$\gl \geq Cd^{-1/4}\log^{3/4}d$ for some large constant $C$.
Specifically, write ${\cal E}$ and ${\cal O}$ for the sets of even
and odd vertices of $\Z^d$ (defined in the obvious way: a vertex
of $\Z^d$ is {\em even} if the sum of its coordinates is even).
Set
$$
\Lambda_M=\Lambda_M^d =[-L,L]^d,~~~
\partial \Lambda_M =[-L,L]^d \setminus [-(L-1),L-1]^d.
$$
For $\lambda>0$, choose ${\bf I}$ from ${\cal I}(\Lambda_M)$ with
$\Pr({\bf I}=I) \propto \lambda^{|I|}$. The main result of
\cite{GK} is that there is a constant $C$ such that if $\lambda
\geq Cd^{-1/4}\log^{3/4}d$ then
$$
\lim_{M\rightarrow\infty}\Pr(\underline{0}\in{\bf I}|{\bf
I}\supseteq
\partial \Lambda_M\cap {\cal E})~>
\lim_{M\rightarrow\infty}\Pr(\underline{0}\in{\bf I}| {\bf
I}\supseteq \partial \Lambda_M\cap {\cal O}).
$$
Thus, roughly speaking, the influence of the boundary on behavior
at the origin persists as the boundary recedes. Informally, this
suggests that for $\gl$ in this range, the typical independent set
chosen according to the hard-core measure is either predominantly
odd or predominantly even. Thus there is a highly unlikely
``bottleneck'' set of balanced independent sets separating the
predominantly odd sets from the predominantly even ones. It is the
existence of this bottleneck that should cause the conductance of
the Glauber dynamics chain to be small (see Section
\ref{sec-proof.of.mn.thm}), and thus cause its mixing time to be
large.

Our main result (Theorem \ref{thm-main}) provides some support for
this belief, verifying it in the case $L=2$; unfortunately,
because of the weak isoperimetry of the torus we cannot hope to
use Theorem \ref{thm-main} to deal with general $L$. (See Remark
\ref{rem-app.to.zd} for further discussion of these issues.)

\medskip

Before stating Theorem \ref{thm-main}, we establish some notation.
From now on, $\gS=(V,E)$ will be a $d$-regular, bipartite graph
with partition classes $\cE$ and $\cO$. Set $N=|V|$ and
$M=|\cE|=|\cO|~(=N/2)$.

For $u, v \in V$ we write $u \sim v$ if there is an edge in $\gS$
joining $u$ and $v$. Set $N(u)=\{w\in V:w \sim u\}$ ($N(u)$ is the
{\em neighbourhood} of $u$) and $N(A)=\cup_{w \in A} N(w)$. For $A
\subseteq \cE$ (or $\cO$) set
$$
[A]=\{x \in V(\gS) : N(x) \subseteq N(A)\};
$$
we think of $[A]$ as an ``external closure'' of $A$. Note that
while $A$ determines $N(A)$, $N(A)$ determines only $[A]$. For
this reason, we find it more convenient at some points in the
sequel to deal with $[A]$ rather than with $A$ itself. Say that
$A$ is {\em small} if $|[A]| \leq M/2$. Define the {\em bipartite
expansion constant} of $\gS$ by
$$
\gd(\gS)=\min \left\{\frac{|N(A)|-|[A]|}{|N(A)|}~:~A \subseteq \cE
\mbox{ small or } A \subseteq \cO \mbox{ small}, A \neq
\emptyset\right\}.
$$
Note that $0 \leq \gd < 1$. (The second inequality is obvious. The
first follows from regularity, which implies that $\gS$ has a
perfect matching, which in turn implies that for all $A \subseteq
\cE$ (or $\cO$), $|A| \leq |N(A)|$.)

All implied constants in $O$ and $\Omega$ notation are independent
of $d$. We use ``$\log$'' throughout for $\log_2$ and ``$\ln$''
for $\log_e$. We write $\exp_2 x$ for $2^x$. We always assume that
$d$ is sufficiently large to support our assertions.

Set
\begin{equation} \label{eq-alpha}
\ga(\gl) = \frac{\log(1+\gl)}{44\left(1+\log
(1+\gl)\right)\log\left(2 + \frac{1}{\log(1+\gl)}\right)}
\end{equation}
and
$$
\gb(\gl) = \frac{\log^2(1+\gl)}{\log (1+\gl) + \log (d^5/\gd)}.
$$
We note for future reference that
\begin{equation} \label{ing-lower.bd.on.alpha}
\gl \geq \frac{1}{\sqrt{d}} ~~\mbox{implies}~~\frac{1}{44}
>\ga(\gl)=\Omega\left(\frac{1}{\sqrt{d}\log d}\right).
\end{equation}
Our main result is

\begin{thm} \label{thm-main}
Let $\gS$ be a $d$-regular, bipartite graph with $N \geq d^2$
vertices and bipartite expansion constant $\gd$. There is a
constant $c > 0$ such that whenever $\gl$ and $\gd$ satisfy
\begin{equation} \label{inq-gen.bd.on.lambda.1}
\gb(\gl) \geq c\max\left\{\frac{\log
(d^5/\gd)}{\sqrt{d}},\frac{\log^2d}{\gd d}\right\}
\end{equation}
we have
$$
\tau_{\cM_\gl(\gS)} \geq \exp_2\left\{\Omega(N
\ga(\gl)\gb(\gl)\gd)\right\}.
$$
\end{thm}

\begin{remark}If we add as an additional hypothesis to Theorem
\ref{thm-main} that $\gS$ has {\em bounded codegree} (that is,
there is a constant $\gk$ independent of $d$ such that each pair
of vertices in $\gS$ has at most $\gk$ common neighbours), then we
can slightly improve our bound on $\gl$ to
\begin{equation} \label{inq-gen.bd.on.lambda.2}
\gb(\gl) \geq c\max\left\{\frac{\log
(d^5/\gd)}{\sqrt{d}},\frac{\log d}{\gd d^2}\right\}.
\end{equation}
We do not present the more complicated argument here.
\end{remark}

Note that since $\gd < 1$, we cannot possibly satisfy
(\ref{inq-gen.bd.on.lambda.1}) for $\gl \leq 1/\sqrt{d}$, so we
may (and will) assume from here on that $\gl \geq 1/\sqrt{d}$.

\medskip

A slightly stronger condition that implies
(\ref{inq-gen.bd.on.lambda.1}) is
\begin{equation} \label{inq-gen.bd.on.lambda.1'}
\log (1 + \gl) \geq c'\max\left\{\frac{\log
(d^5/\gd)}{d^{1/4}},\frac{\log d \sqrt{\log (d^5/\gd)}}{\sqrt{\gd
d}}\right\}
\end{equation}
where the constant $c'$ depends on $c$, from which we can more
clearly see the tradeoff between $\gl$ and $\gd$. From
(\ref{inq-gen.bd.on.lambda.1'}) we may also read off the following
corollary of Theorem \ref{thm-main} addressing Glauber dynamics
for sampling a uniform independent set ($\gl =1$).

\begin{cor}
Let $\gS$ satisfy the conditions of Theorem \ref{thm-main}. There
is a constant $c > 0$ such that whenever $\gd \geq c\log^3 d/d$ we
have
$$
\tau_{\cM_1(\gS)} \geq \exp_2\left\{\Omega\left(\frac{N\gd}{\log
d}\right)\right\}.
$$
\end{cor}

\medskip

As an application of Theorem \ref{thm-main}, we consider the case
$\gS=Q_d$, the $d$-dimensional hypercube. This is the $d$-regular
bipartite graph on vertex set $\{0,1\}^d$ in which two vertices
are adjacent if they differ on exactly one coordinate.
%
%
%
%
The hypercube satisfies $\gd(Q_d)\geq \Omega(1/\sqrt{d})$ (see,
e.g. \cite[Lemma 1.3]{K-S}; this bound can also be derived from
isoperimetric inequalities of Bezrukov \cite{Bezrukov} and
K\"orner and Wei \cite{KornerWei}) and so
if $c'>0$ is a suitably large constant (depending on the constant
$c$ provided by Theorem \ref{thm-main}) then
(\ref{inq-gen.bd.on.lambda.1}) is satisfied as long as $\gl \geq
c'd^{-1/4}\log^{3/2} d$. So the following is a corollary of
Theorem \ref{thm-main}.

\begin{cor} \label{cor-hypercube}
There are constants $c, c' >0$ such that whenever $\gl \geq
cd^{-1/4}\log^{3/2} d$ we have
$$
\tau_{\cM_\gl(Q_d)} \geq
\exp_2\left\{\frac{c'2^d\log^3(1+\gl)}{\sqrt{d}\left(1+\log
(1+\gl)\right)\left(c_3\log d+\log (1+\gl)\right)\log\left(2 +
\frac{1}{\log(1+\gl)}\right)}\right\}.
$$
In particular,
$$
\tau_{\cM_\gl(Q_d)} \geq \left\{\begin{array}{ll}
\exp_2\left\{\frac{2^d\log^3(1+\gl)}{\sqrt{d}\log^2d}\right\} &
\mbox{if ~$c d^{-1/4} \log^{3/2} d \leq \gl \leq
O(1)$}, \\
\\
\exp_2\left\{\frac{2^d\log^2(1+\gl)}{\sqrt{d}\log d}\right\} &
\mbox{if ~$\Omega(1) \leq \gl \leq
O(d)$}, \\
\\
\exp_2\left\{\frac{2^d\log(1+\gl)}{\sqrt{d}}\right\} & \mbox{if
~$\Omega(d) \leq \gl$}.
                                  \end{array}
                            \right.
$$
\end{cor}

\begin{remark}
Using (\ref{inq-gen.bd.on.lambda.2}) in place of
(\ref{inq-gen.bd.on.lambda.1}) (which we may do, since $Q_d$ has
bounded codegree) we may improve the bound on $\gl$ in Corollary
\ref{cor-hypercube} to $\gl \geq cd^{-1/4} \log d$.
\end{remark}

\begin{remark} \label{rem-app.to.zd}
Let us return to $T_{L,d}$, the even discrete torus. Since $Q_d$
is easily seen to be isomorphic to $T_{2,d}$, Corollary
\ref{cor-hypercube} gives an exponential lower bound on
$\tau_{\cM_\gl(T_{2,d})}$ for sufficiently large $d$ whenever $\gl
\geq \omega(d^{-1/4}\log^{3/2} d)$. Unfortunately, the best bound
we can obtain on the bipartite expansion constant of $T_{L,d}$ is
$\gd(T_{L,d})\geq \Omega(1/Ld)$ (see, e.g \cite{GK}), so we cannot
use Theorem \ref{thm-main} to obtain {\em any} lower bound on
$\gl$ independent of $L$ beyond which $\tau_{\cM_\gl(T_{L,d})}$ is
large for all even $L\geq 4$ and sufficiently large $d$. However,
subsequent to the completion of this paper, a strategy specific to
the torus has been employed in \cite{Galvin} to show that for all
even $L\geq 4$, $\gl \geq \omega(d^{-1/4}\log^{3/4}d)$ and $d$
sufficiently large,
$$
\tau_{\cM_\gl(T_{L,d})} \geq \exp_2\left\{\frac{L^{d-1}}{d^3\log
^2 L}\right\}.
$$
\end{remark}

\section{Proof of Theorem \ref{thm-main}} \label{sec-proof.of.mn.thm}

The notion of conductance, introduced in \cite{JerrumSinclair},
can be used to analyze the behavior of $\tau_{\cM_\gl}$. Let $\cM$
be a Markov chain on state space $\gO$ with transition matrix $P$
and stationary distribution $\pi$. For $\go_1, \go_2 \in \gO$ and
$A,B \subseteq \gO$, set
$$
Q(\go_1,\go_2)=\pi(\go_1)P(\go_1,\go_2) ~~~~\mbox{and}~~~~
Q(A,B)=\sum_{\go_1\in A, ~\go_2 \in B} Q(\go_1,\go_2).
$$
For $\emptyset \neq S \subseteq \gO$, define the conductance of
$S$ as
$$
\Phi(S)=\frac{Q(S,\gO \setminus S)}{\pi(S)}.
$$
We may interpret $\Phi(S)$ as the probability under $\pi$ that the
chain escapes from $S$ in one step, given that it is in $S$.
Define the conductance of $\cM$ as
$$
\Phi_{\cM}=\min_{0 < \pi(S) \leq \frac{1}{2}} \Phi(S).
$$
We may then bound the mixing time of $\cM$ by
\begin{equation} \label{inq-mixing.vs.conductance}
\tau_{\cM} \geq
\left(\frac{1}{2}-\frac{1}{e}\right)\frac{1}{\Phi_{\cM}}
\end{equation}
(see e.g. \cite{DyerFriezeJerrum}, where the above bound is
derived without assuming time-reversibility of the chain $\cM$).
Thus to show that the mixing time is large, it is enough to
exhibit a single $S$ with small conductance.

Throughout this section we fix $\gS$ satisfying the conditions of
Theorem \ref{thm-main}. Set
$$
\cI_\cE = \{I \in \cI(\gS)~:~|I\cap \cE| > |I\cap \cO|\},
$$
define $\cI_\cO$ analogously, and set $\cI_b = \cI(\gS) \setminus
(\cI_\cE \cup \cI_\cO)$ ($\cI_b$ is the set of {\em balanced}
independent sets). Without loss of generality, assume
$\pi_\gl(\cI_\cE) \leq 1/2$. Because Glauber dynamics changes the
size of an independent set by at most one at each step, we have
that if $I \in \cI_\cE, J \not \in \cI_\cE$ satisfy $P_\gl(I,J)
\neq 0$, then $J \in \cI_b$. It follows that
\begin{eqnarray}
Q(\cI_\cE, \Omega \setminus \cI_\cE) & = & \sum_{I \in \cI_\cE, J
\not \in
\cI_\cE} \pi_\gl(I)P_\gl(I,J) \nonumber \\
& = & \sum_{I \in \cI_\cE, J \not \in \cI_\cE}
\pi_\gl(J)P_\gl(J,I) \label{using.reversibility} \\
& = & \sum_{I \in \cI_\cE, J \in \cI_b}
\pi_\gl(J)P_\gl(J,I) \nonumber \\
& \leq & \pi_\gl(\cI_b). \nonumber
\end{eqnarray}
The simplest way to see (\ref{using.reversibility}) is to use the
fact that $\cM_\gl$ is time-reversible (that is, that
$\pi_\gl(I)P_\gl(I,J)=\pi_\gl(J)P_\gl(J,I)$ for all $I, J \in
\cI$); but note that more generally if $\cM$ is a (not necessarily
time-reversible) Markov chain on finite state space $\Omega$ with
transition matrix $P$ and stationary distribution $\pi$ then
$$
\sum_{\go_1 \in S,\go_2 \not \in S}
\pi(\go_1)P(\go_1,\go_2)=\sum_{\go_1 \in S,\go_2 \not \in S}
\pi(\go_2)P(\go_2,\go_1)
$$
for all $S \subseteq \Omega$. Now using the trivial lower bound
$w_\gl(\cI_\cE) \geq (1+\gl)^M$ (recall that for $\cJ \subseteq
\cI$, $\omega_\gl(\cJ)=\sum_{J \in \cJ}\gl^{|J|}$) we obtain
\begin{equation} \label{cond.bound}
\Phi_{\cM_\gl} \leq \Phi(\cI_\cE) \leq
\frac{\pi_\gl(\cI_b)}{\pi_\gl(\cI_\cE)} =
\frac{w_\gl(\cI_b)}{w_\gl(\cI_\cE)} \leq
\frac{w_\gl(\cI_b)}{(1+\gl)^M}.
\end{equation}
Thus (recalling
(\ref{inq-mixing.vs.conductance})) to show that $\tau_{\cM_\gl}$
is large, it is enough to show that $w_\gl(\cI_b)$ is small. We
may think of $\cI_b$ as a ``bottleneck'' set through which any run
of the chain must pass in order to mix; if the bottleneck has low
measure, the mixing time is high.

We will actually consider a larger ``bottleneck'' set. Set
$$
\cI^{triv} = \{I \in \cI ~:~ |I\cap \cE|,|I \cap \cO| \leq
\ga(\gl)M\}
$$
and
$$
\cI^{nt} = \{I \in \cI ~:~ \min\{|I\cap \cE|,|I \cap \cO|\}\geq
\ga(\gl)M\},
$$
where $\ga(\gl)$ is as defined in (\ref{eq-alpha}). Note that
$\cI_b \subseteq \cI^{triv} \cup \cI^{nt}$. We will show that as
long as $\gl$ satisfies (\ref{inq-gen.bd.on.lambda.1}),
\begin{equation} \label{inq-bounding.triv.plus.nt}
w_\gl(\cI^{triv} \cup \cI^{nt}) \leq (1+\gl)^M
\exp_2\left\{-\Omega\left(M\ga(\gl)\gb(\gl)\gd\right)\right\},
\end{equation}
from which Theorem \ref{thm-main} follows via
(\ref{inq-mixing.vs.conductance}) and (\ref{cond.bound}).

\medskip

Dealing with $w_\gl(\cI^{triv})$ is relatively straightforward. We
begin by observing that
\begin{equation} \label{inq-alphalogalpha}
4\ga(\gl)\log \frac{1}{\ga(\gl)} \leq
\frac{\log(1+\gl)}{2\left(1+\log(1+\gl)\right)}.
\end{equation}
To see this, first set
$$
\grg(\gl) = \frac{\log(1+\gl)}{1+\log(1+\gl)}.
$$
Note that for all $\gl > 0$, $0 < \grg(\gl) < 1$. We have
$\ga(\gl)=\grg(\gl)/(44\log(1+1/\grg(\gl)))$ and so
(\ref{inq-alphalogalpha}) is equivalent to
$$
\frac{\grg(\gl)}{11\log\left(1+\frac{1}{\grg(\gl)}\right)}
\log\left(\frac{44\log\left(1+\frac{1}{\grg(\gl)}\right)}{\grg(\gl)}\right)
\leq \frac{\grg(\gl)}{2}
$$
which is in turn equivalent to
$$
44\log \left(1+\frac{1}{\grg(\gl)}\right) \leq
\grg(\gl)\left(1+\frac{1}{\grg(\gl)}\right)^{11/2}.
$$
That this inequality holds for all $0 < \grg(\gl) < 1$ is a
routine calculus exercise. Note also that for $0 < x < 1/e$,
\begin{equation} \label{inq-entropy}
x \leq H(x) \leq 2x\log \frac{1}{x}
\end{equation}
(where recall $H(x)=-x\log x -(1-x)\log(1-x)$
is the usual binary entropy function).
Finally, we use a result concerning the sums of binomial
coefficients which follows from the Chernoff bounds \cite{Ch} (see
also \cite{RG}, p.11):
\begin{equation} \label{inq-binomial}
\sum_{i=0}^{[cN]}{N\choose i} \leq 2^{H(c)N}~~~~~\mbox{for $c \leq
\frac{1}{2}$},
\end{equation}
where $[x]$ denotes the integer part of $x$.

Now
with the inequalities justified below, we have
\begin{eqnarray}
w_\gl(\cI^{triv}) & \leq &
{M\choose \leq \ga(\gl)M}^2(1+\gl)^{2\ga(\gl)M} \nonumber \\
& \leq & \exp_2\left\{2H(\ga(\gl))M+2\ga(\gl)M\log(1+\gl)\right\}
\label{using-bin.entropy} \\
& \leq & \exp_2\left\{2MH(\ga(\gl))(1+\log(1+\gl))\right\}
\label{using-entropy.lb} \\
& \leq &
\exp_2\left\{4M\ga(\gl)\log\frac{1}{\ga(\gl)}(1+\log(1+\gl))\right\}
\label{using-entropy.ub} \\
& \leq & \exp_2\left\{M\frac{\log(1+\gl)}{2}\right\}
\label{using-alphalogalpha} \\
& \leq & (1+\gl)^{\frac{M}{2}} \label{inq-bounding.triv}.
\end{eqnarray}
Here (and throughout) we use ${n \choose \leq k}$ for $\sum_{i
\leq k}{n \choose i}$. In (\ref{using-bin.entropy}), we are using
(\ref{inq-binomial}), which is applicable by
(\ref{ing-lower.bd.on.alpha}). In (\ref{using-entropy.lb}) we are
using the first inequality in (\ref{inq-entropy}) and in
(\ref{using-entropy.ub}) we are using the second (again, both of
these are applicable by (\ref{ing-lower.bd.on.alpha}).) Finally in
(\ref{using-alphalogalpha}) we are using
(\ref{inq-alphalogalpha}).

Bounding $w_\gl(\cI^{nt})$ requires much more work. We begin by
enlarging $\cI^{nt}$ slightly. Say that $I \in \cI(\gS)$ is {\em
small on $\cE$} if $|[I\cap \cE]| \leq M/2$ (recall that for $A
\subseteq \cE$, $[A]=\{v \in \cO:N(v)\subseteq A\}$), and set
$$
\cI^{nt}_{\cE} =\{I \in \cI^{nt}: I \mbox{ is small on $\cE$}\}.
$$
Define {\em small on $\cO$} and $\cI^{nt}_{\cO}$ similarly. A
simple argument, based on the fact that $\gS$ has a perfect
matching, shows that any $I \in \cI(\gS)$ must be small on at
least one of $\cE$, $\cO$, and so we have
$$
w_\gl(\cI^{nt}) \leq
2\max\left\{w_\gl(\cI^{nt}_{\cE}),w_\gl(\cI^{nt}_{\cO})\right\}.
$$
We may assume without loss of generality that
$$
w_\gl(\cI^{nt}_{\cE}) =
\max\left\{w_\gl(\cI^{nt}_{\cE}),w_\gl(\cI^{nt}_{\cO})\right\}
$$
so that it is enough to show that
$$
w_{\gl}(\cI^{nt}_{\cE}) \leq (1+\gl)^M \exp_2\left\{-\Omega(M
\ga(\gl) \gb(\gl) \gd)\right\}.
$$

For each $a \geq \ga(\gl)M$ and $g \geq a$ set
$$
\cA(a,g)=\{A \subseteq \cE: |[A]|=a, |N(A)|=g\}
$$
and set
$$
\cI(a,g)=\{I \in \cI^{nt}_{\cE}:I\cap \cE \in \cA(a,g)\}.
$$
We have
\begin{eqnarray*}
w_{\gl}(\cI^{nt}_{\cE}) & \leq &
\sum_{a \geq \ga(\gl)M, ~g \geq a}w_\gl(\cI(a,g)) \\
& \leq & \sum_{a \geq \ga(\gl)M, ~g \geq a}w_\gl(\cA(a,g))(1+\gl)^{M-g}\\
& \leq & (1+\gl)^M\sum_{a \geq \ga(\gl)M, ~g \geq a}
w_\gl(\cA(a,g))(1+\gl)^{-g} \\
& \leq & (1+\gl)^M M^2\max_{a \geq \ga(\gl)M, ~g \geq a}
w_\gl(\cA(a,g))(1+\gl)^{-g}
\end{eqnarray*}

The key now is to upper bound $w_\gl(\cA(a,g))$. The following
theorem (whose proof is given in Section
\ref{sec-proof.of.mn.lem}) is based on ideas of A. Sapozhenko
\cite{Sapozhenko,Sap2}.

\begin{thm} \label{lem-main}
Let $\gS$ be any graph satisfying the assumptions of Theorem
\ref{thm-main}. We have
$$
w_\gl(\cA(a,g)) \leq (1+\gl)^g
\exp_2\left\{-\Omega\left((g-a)\gb(\gl)\right)\right\}.
$$
for any $a \geq \ga(\gl)M$ and any $\gl$ satisfying
(\ref{inq-gen.bd.on.lambda.1}).
\end{thm}

For $a \geq \ga(\gl)M$ and $g \geq a$ we have $g-a \geq \gd g \geq
M\ga(\gl) \gd$ and so
\begin{eqnarray}
w_{\gl}(\cI^{nt}_{\cE}) & \leq &
(1+\gl)^MM^2\exp_2\left\{-\Omega(M\ga(\gl)\gb(\gl)\gd)\right\}
\nonumber \\
& \leq & (1+\gl)^M
\exp_2\left\{-\Omega(M\ga(\gl)\gb(\gl)\gd)\right\}.
\label{inq-bounding.nt}
\end{eqnarray}
To see that the factor of $M^2$ may be absorbed into the exponent,
note that by hypothesis, $2M\geq d^2$ and so $M^2 \leq
\exp_2\left\{O(M\log d/d^2)\right\}$, and that combining
(\ref{ing-lower.bd.on.alpha}) and (\ref{inq-gen.bd.on.lambda.1})
we have $\ga(\gl)\gb(\gl)\gd \geq \Omega(d^{-3/2}\log d)$.

Combining (\ref{inq-bounding.nt}) and (\ref{inq-bounding.triv}) we
get (\ref{inq-bounding.triv.plus.nt}) and hence Theorem
\ref{thm-main}.

\section{Proof of Theorem \ref{lem-main}}
\label{sec-proof.of.mn.lem}

For $u, v \in V$ and $A,B \subseteq V$ we write $\nabla(A)$ for
the set of edges having one end in $A$ and (if $A \cap B =
\emptyset$) $\nabla(A,B)$ for the set of edges having one end in
each of $A, B$. We also write $d_A(v)$ for $|N(v) \cap A|$.

Throughout this section, we fix $\gS$ satisfying the assumptions
of Theorem \ref{thm-main}. We also fix $a$ and $g$, but we do not
assume $a \geq \ga(\gl)M$. We write $\cA$ for $\cA(a,g)$. Given $A
\in \cA$ we always write $G$ for $N(A)$ and set $t=g-a$. Note that
for $A \in \cA$,
\begin{equation} \label{eq-nabla}
|\nabla(G, \cE \setminus [A])| =dg-da=td.
\end{equation}

The proof of Theorem \ref{lem-main} involves the idea of
``approximation''. We begin with an informal outline. To bound
$w_\gl(\cA)$, we produce a small set $\cU$ with the properties
that each $A \in \cA$ is ``approximated'' (in an appropriate
sense) by some $U \in \cU$, and for each $U \in {\cal U}$, the
total weight of those $A \in \cA$ that could possibly be
``approximated'' by $U$ is small. (Each $U \in {\cal U}$ will
consist of two parts; one each approximating $G$ and $A$.) The
product of the bound on $|\cU|$ and the bound on the weight of
those $A \in \cA$ that may be approximated by any $U$ is then a
bound on $w_\gl(\cA)$. The set $\cU$ is itself produced by an
approximation process --- we first produce a small set $\cV$ with
the property that each $A \in \cA$ is ``weakly approximated'' (in
an appropriate sense) by some $V \in \cV$, and then show that for
each $V$ there is a small set $\cW(V)$ with the property that for
each $A \in \cA$ that is ``weakly approximated'' by $V$, there is
a $W \in \cW(V)$ which approximates $A$; we then take $\cU=\cup_{V
\in \cV} \cW(V)$. (Each $V \in \cV$ will consist of a single
part.)

The main inspiration for the proof of Theorem \ref{lem-main} is
the work of A. Sapozhenko, who, in \cite{Sap2}, gave a relatively
simple derivation for the asymptotics of the number of independent
sets in $Q_d$ (in the notation of (\ref{eq-hcmeasure}), this is
the asymptotics of $Z_{\gl}(Q_d)$ with $\gl =1$), earlier derived
in a more involved way in \cite{K-S}. Our Lemma
\ref{lem-psi.approx} is a modification of a lemma in
\cite{Sapozhenko}, and our overall approach is similar to
\cite{Sap2}. See e.g. \cite{G} for another recent application of
these ideas.

We now begin the formal discussion of Theorem \ref{lem-main} by
introducing the two notions of approximation that we will use,
beginning with the weaker notion. A {\em covering approximation}
for $A \subseteq \cE$ is a set $F_0 \in 2^\cO$ satisfying
$$
F_0 \subseteq G,~N(F_0) \supseteq [A].
$$

The second notion of approximation depends on a parameter $\psi$,
$1 \leq \psi \leq d/2$. A {\em $\psi$-approximation} for $A
\subseteq \cE$ is a pair $(F,S) \in 2^\cO \times 2^\cE$ satisfying
\begin{equation} \label{cond-psiapprox1}
F \subseteq G,~~S \supseteq [A],
\end{equation}
\begin{equation} \label{cond-psiapprox2}
d_F(u) \geq d-\psi~~~\forall u \in S
\end{equation}
and
\begin{equation} \label{cond-psiapprox3}
d_{\cE \setminus S}(v) \geq d-\psi~~~\forall v \in \cO \setminus
F.
\end{equation}
Note that if $x\in [A]$ then $N(x) \subseteq G$, and if $y \in \cO
\setminus G$ then $N(y) \subseteq \cE \setminus [A]$. If we think
of $S$ as ``approximate $[A]$'' and $F$ as ``approximate $G$'',
(\ref{cond-psiapprox2}) says that if $x \in \cE$ is in
``approximate $[A]$'' then almost all of its neighbours are in
``approximate $G$'', while (\ref{cond-psiapprox3}) says that if
$y\in \cO$ is not in ``approximate $G$'' then almost all of its
neighbours are not in ``approximate $[A]$''.

Before continuing, we note a property of $\psi$-approximations
that will be of use later.
\begin{lemma} \label{lem-s.small}
If $(F,S)$ is a $\psi$-approximation for $A \in \cA$ then
\begin{equation} \label{inq-bounding.s.by.f}
|S| \leq |F| + \frac{2t\psi}{d-\psi}.
\end{equation}
\end{lemma}

\noindent {\em Proof: }Observe that $|\nabla(S,G)|$ is bounded
above by $d|F| + \psi|G \setminus F|$ and below by $d|[A]| +
(d-\psi)|S \setminus [A]| = d|S| - \psi|S \setminus [A]|$, giving
$$
|S| \leq |F| + \frac{\psi|(G \setminus F) \cup (S \setminus
[A])|}{d},
$$
and that each $u \in (G \setminus F) \cup (S \setminus [A])$
contributes at least $d-\psi$ edges to $\nabla(G,\cE\setminus
[A])$, a set of size $td$, giving
$$
|(G \setminus F) \cup (S \setminus [A])| \leq \frac{2td}{d-\psi}.
$$
These two observations together give (\ref{inq-bounding.s.by.f}).
\qed

\medskip

There are three parts to the proof of Theorem \ref{lem-main}.
Lemma \ref{lem-phi.approx} is the first ``approximation'' step,
producing a small family $\cV$ of covering approximations for
$\cA$. Lemma \ref{lem-psi.approx} is the second ``approximation''
step, refining the covering approximations to produce a family
$\cW$ of $\psi$-approximations for $\cA$. Finally, Lemma
\ref{lem-reconstruction} is the ``reconstruction'' step, bounding
the weight of the set of $A$'s that could possibly be
$\psi$-approximated by a member of $\cW$. We now state the three
relevant lemmas. We will then derive Theorem \ref{lem-main} before
turning to the proofs of the approximation and reconstruction
lemmas.

\begin{lemma} \label{lem-phi.approx}
There is a $\cV = \cV(a,g) \subseteq 2^\cO$ with
$$
|\cV| \leq {M \choose \leq \frac{2g\log d}{d}}
$$
such that each $A \in \cA$ has a covering approximation in $\cV$.
\end{lemma}

\begin{lemma} \label{lem-psi.approx}
For any $F_0 \in \cV$ and $1 \leq \psi \leq d/2$ there is a
$\cW=\cW(F_0,\psi, a, g) \subseteq 2^\cO \times 2^\cE$ with
$$
|\cW| \leq {2g\log d \choose \leq \frac{2g}{d}}{2d^3g\log d
\choose \leq \frac{2t}{\psi}}{2g\log d \choose \leq
\frac{td}{(d-\psi)\psi}}
$$
such that any $A \in \cA$ for which $F_0$ is a covering
approximation has a $\psi$-approximation in $\cW$.
\end{lemma}

\begin{lemma} \label{lem-reconstruction}
Given $1 \leq \psi \leq d/2$ and $1 \geq \grg >
\frac{-2\psi}{d-\psi}$, for each $(F,S) \in 2^\cO \times 2^\cE$
that satisfies (\ref{inq-bounding.s.by.f}) we have
\begin{equation} \label{inq-cost.of.reconstruction}
\sum w_\gl(A) \leq \max \left\{(1+\gl)^{g-\grg t}, {3dg \choose
\leq \frac{2t\psi}{d-\psi}+\grg t}(1+\gl)^{g-t}\right\}
\end{equation}
where the sum is over all those $A$'s in $\cA$ satisfying $F
\subseteq G$ and $S \supseteq [A]$.
\end{lemma}

\bigskip

Before turning to the proofs of Lemmas \ref{lem-phi.approx},
\ref{lem-psi.approx} and \ref{lem-reconstruction}, we use them to
obtain Theorem \ref{lem-main}. Throughout, we will use (usually
without comment) a simple observation about sums of binomial
coefficients: if $k=o(n)$, we have
\begin{eqnarray}
\sum_{i \leq k} {n \choose i} & \leq & (1+O(k/n)){n \choose k} \nonumber \\
& \leq & (1+O(k/n))(en/k)^k \nonumber \\
& \leq & \exp_2\left\{(1+o(1))k\log(n/k)\right\}. \label{binomial}
\end{eqnarray}

Take $\psi=\sqrt{d}$ and
$$
\grg = \frac{\log(1+\gl)-\frac{\sqrt{d}}{d-\sqrt{d}}\log
(d^5/\gd)}{\log(1+\gl)+\log (d^5/\gd)}.
$$
Note that for this choice of $\psi$ and $\grg$ we have $\grg >
\frac{-\psi}{d-\psi}$, and so
$$
\log \frac{3d}{\gd\left(\frac{2\psi}{d-\psi}+\grg\right)} \leq
\frac{1}{2}\log (d^5/\gd).
$$
The bound in (\ref{inq-cost.of.reconstruction}) is therefore at
most
$$
(1+\gl)^g \exp_2\left\{\max \left\{ - \grg t\log(1+\gl),~
t\left(\frac{2\psi}{d-\psi}+\grg \right)\log (d^5/\gd)
\right\}\right\}.
$$
(Here we have used (\ref{binomial})). For our choice of $\psi$ and
$\grg$ this is at most
$$
(1+\gl)^g
\exp_2\left\{t\frac{\log(1+\gl)\frac{\sqrt{d}}{d-\sqrt{d}}\log
(d^5/\gd) - \log^2(1+\gl)}{\log(1+\gl)+\log (d^5/\gd)} +
t\frac{\sqrt{d}\log (d^5/\gd)}{d-\sqrt{d}}\right\},
$$
which in turn is at most
\begin{equation} \label{inq-cost.of.reconstruction.specific}
(1+\gl)^g \exp_2\left\{O\left(t\frac{\log
(d^5/\gd)}{\sqrt{d}}\right) -
t\frac{\log^2(1+\gl)}{\log(1+\gl)+\log (d^5/\gd)}\right\}.
\end{equation}
The bounds in Lemmata \ref{lem-psi.approx} and
\ref{lem-phi.approx} are at most
\begin{equation} \label{inq-cost.of.psi.phi.specific}
\exp_2\left\{O\left(g\frac{\log d}{d}+t\frac{\log
(d^5/\gd)}{\sqrt{d}}\right)\right\}
~~~~\mbox{and}~~~~\exp_2\left\{O\left(g\frac{\log^2
d}{d}\right)\right\}
\end{equation}
respectively. For the latter bound, we are using the assumption
$a\geq \ga(M)$ of Theorem \ref{lem-main} and the fact that $\gl \geq
1/\sqrt{d}$, which together imply (via
(\ref{ing-lower.bd.on.alpha})) that
$$
\frac{Md}{g\log d} \leq d^{3/2}.
$$
Combining (\ref{inq-cost.of.psi.phi.specific}) with
(\ref{inq-cost.of.reconstruction.specific}), we get
$$
w_\gl(\cA) \leq (1+\gl)^g \exp_2\left\{O\left(g\frac{\log^2 d}{d}+
t\frac{\log (d^5/\gd)}{\sqrt{d}}\right)
                    -t\frac{\log^2(1+\gl)}{\log(1+\gl)+\log
                    (d^5/\gd)}\right\}.
$$
Noting that $t \geq \gd g$ always, we find that if $\gl$ satisfies
(\ref{inq-gen.bd.on.lambda.1}) with a suitably large constant $c$,
then
$$
w_\gl(\cA) \leq (1+\gl)^g
\exp_2\left\{-\Omega\left(t\frac{\log^2(1+\gl)}{\log (1+\gl) +
\log (d^5/\gd)}\right)\right\}.
$$
and so we get Theorem \ref{lem-main}.

\bigskip

We now turn to the proofs of Lemmata \ref{lem-phi.approx},
\ref{lem-psi.approx} and \ref{lem-reconstruction}.

\bigskip

\noindent {\em Proof of Lemma \ref{lem-phi.approx}: } We appeal to
a special case of a fundamental result due to Lov\'asz
\cite{Lovasz} and Stein \cite{Stein}. For a bipartite graph $\gG$
with bipartition $P \cup Q$, we say that $Q^\prime \subseteq Q$
{\em covers} $P$ if each $p \in P$ has a neighbour in $Q^\prime$.

\begin{lemma} \label{lovaszstein}
If $\gG$ as above satisfies $|N(x)| \geq p$ for each $x \in P$ and
$|N(y)| \leq q$ for each $y \in Q$, then $P$ is covered by some
$Q^\prime \subseteq Q$ with
$$|Q^\prime| \leq (|Q|/p)(1 + \ln q).$$
\end{lemma}

Applying the lemma with $\gG$ the subgraph of $\gS$ induced by
$[A]\cup G$, $P=[A]$, $Q=G$ and $p=q=d$, we find that each $A \in
\cA$ has a covering approximation of size at most $2g\log d/d$.
Taking $\cV$ to be the set of all subsets of $\cO$ of size at most
$2g\log d/d$, the lemma follows. \qed

\bigskip

\noindent {\em Proof of Lemma \ref{lem-psi.approx}: } We describe
an algorithm, which we refer to as the {\em degree algorithm},
which produces for input $(F_0,A) \in 2^\cO \times 2^\cE$ for
which $F_0$ is a covering approximation of $A$ (i.e., with
$N(F_0)\supseteq [A]$), an output $(F,S) \in 2^\cO \times 2^\cE$
which is a $\psi$-approximation for $A$ (i.e, which satisfies
(\ref{cond-psiapprox1}), (\ref{cond-psiapprox2}) and
(\ref{cond-psiapprox3})). The idea for the algorithm is from
\cite{Sapozhenko}. To begin, fix a linear ordering $\ll$ of $V$.

\medskip

\noindent {\bf Step $1$: }If $\{u \in [A]: d_{G \setminus F_0}(u)
>d/2 \}
\neq \emptyset$, pick the smallest (with respect to $\ll$) $u$ in
this set and update $F_0$ by $F_0 \longleftarrow F_0 \cup N(u)$.
Repeat this until $\{u \in [A]: d_{G \setminus F_0}(u) > d/2\} =
\emptyset$. Then set $F_1=F_0$ and $S_1=\{u \in \cE:d_{F_1}(u)
\geq d-d/2\}$ and go to Step $2$.

\medskip

\noindent {\bf Step $2$: }If $\{v \in \cO \setminus G: d_{S_1}(v)
>\psi \} \neq \emptyset$, pick the smallest (with respect to
$\ll$) $v$ in this set and update $S_1$ by $S_1 \longleftarrow S_1
\setminus N(v)$. Repeat this until $\{v \in \cO \setminus G:
d_{S_1}(v)
>\psi \} = \emptyset$. Then set $S_2=S_1$ and
$F_2= \{v \in \cO :d_{S_2}(v) > \psi \}$ and go to Step $3$.

\medskip

\noindent {\bf Step $3$: }If $\{w \in [A]: d_{G \setminus F_2}(w)
>\psi \}
\neq \emptyset$, pick the smallest (with respect to $\ll$) $w$ in
this set and update $F_2$ by $F_2 \longleftarrow F_2 \cup N(w)$.
Repeat this until $\{w \in [A]: d_{G \setminus F_2}(w) > \psi \} =
\emptyset$. Then set $F=F_2$ and $S=S_2 \cap \{w \in \cE:d_{F}(w)
\geq d-\psi\}$ and stop.

\begin{claim} \label{algoanalysisoutput}
The output of the degree algorithm is a $\psi$-approximation for
$A$.
\end{claim}

\noindent {\em Proof: }To see that $F \subseteq G$ and $S
\supseteq [A]$, first observe that $S_1 \supseteq [A]$ (or Step
$1$ would not have terminated). We then have $S_2 \supseteq [A]$
(since Step $2$ deletes from $S_1$ only neighbours of $\cO
\setminus G$), and $F_2 \subseteq G$ (or Step $2$ would not have
terminated). Finally, $F\subseteq G$ (since the vertices added to
$F_2$ in Step $3$ are all in $G$) and $S \supseteq [A]$ (or Step
$3$ would not have terminated)

By the definition of $S$, (\ref{cond-psiapprox2}) is satisfied. To
verify (\ref{cond-psiapprox3}), note that by definition of $F_2$,
if $y \in \cO \setminus F_2$ then $d_{\cE \setminus S_2}(y)\geq d
-\psi$. That $y \in \cO \setminus F$ implies $d_{\cE \setminus
S}(y)\geq d -\psi$ now follows from the fact that $F_2 \subseteq
F$ and $S_2 \supseteq S$. \qed

\medskip

\begin{remark}
The alert reader may have noticed that if we replace $d/2$ by
$\psi$ in Step $1$, then the output of Step $2$ is already a
$\psi$-approximation for $A$. The three-step algorithm, however,
is needed to obtain the right bound on $\gb(\gl)$ in Theorem
\ref{thm-main}; see Remark \ref{rem-three.step} following the
proof of Claim \ref{algoanalysisnum}.
\end{remark}

\begin{claim} \label{algoanalysisnum}
Fix $F_0 \in \cV$. The degree algorithm has at most
$$
{2g\log d \choose \leq \frac{2g}{d}}{2d^3g\log d \choose \leq
\frac{2t}{\psi}}{2g\log d \choose \leq \frac{td}{(d-\psi)\psi}}
$$
outputs as the input runs over those $(F_0, A)$ for which $A \in
\cA$ and $F_0$ is a covering approximation for $A$.
\end{claim}

\noindent Taking ${\cal W}$ to be the set of all possible outputs
of the algorithm, the lemma follows.

\medskip

\noindent {\em Proof of Claim \ref{algoanalysisnum}: } The output
of the algorithm is determined by the set of $u$'s whose
neighbourhoods are added to $F_0$ in Step $1$, the set of $v$'s
whose neighbourhoods are removed from $S_1$ in Step $2$, and the
set of $w$'s whose neighbourhoods are added to $F_2$ in Step $3$.

Each iteration in Step $1$ removes at least $d/2$ vertices from
$G\setminus F$, a set of size at most $g$, so there are at most
$2g/d$ iterations. The $u$'s in Step $1$ are all drawn from $[A]$
and hence $N(F_0)$, a set of size at most $d|F_0|\leq 2g\log d$.
So the total number of outputs for Step $1$ is at most
\begin{equation} \label{inq-cost.of.step.1}
{2g\log d \choose \leq \frac{2g}{d}}.
\end{equation}

At the start of Step $2$, each $x \in S_1 \setminus [A]$
contributes at least $d/2$ edges to $\nabla(G, \cE \setminus
[A])$, by (\ref{eq-nabla}) a set of size $dt$, so $|S_1\setminus
[A]|\leq 2t$. Each $v$ used in Step $2$ reduces this by at least
$\psi$, so there are at most $2t/\psi$ iterations. Each $v$ is
drawn from $N(S_1)$, a set which is contained in the fourth
neighbourhood of $F_0$ ($S_1 \subseteq N(G)$ by construction of
$S_1$, $G=N(A)$ and $A \subseteq N(F_0)$) and so has size at most
$d^4|F_0|\leq 2d^3 g\log d$. So the total number of outputs for
Step $2$ is
\begin{equation} \label{inq-cost.of.step.2}
{2d^3g\log d \choose \leq \frac{2t}{\psi}}.
\end{equation}

At the start of Step $3$, each $y \in G \setminus F_2$ contributes
at least $d-\psi$ edges to $\nabla(G, \cE \setminus [A])$, so
$|G\setminus F_2|\leq dt/(d-\psi)$. Each $w$ used in step $3$
reduces this by at least $\psi$, so there are at most
$dt/((d-\psi)\psi)$ iterations. As in Step $1$, the $w$'s are all
drawn from a set of size at most $2g\log d$, so the total number
of outputs for Step $1$ is at most
\begin{equation} \label{inq-cost.of.step.3}
{2g\log d \choose \leq \frac{td}{(d-\psi)\psi}}.
\end{equation}
Combining (\ref{inq-cost.of.step.1}), (\ref{inq-cost.of.step.2})
and (\ref{inq-cost.of.step.3}), the claim follows. \qed

\medskip

\begin{remark} \label{rem-three.step}
The bound in Claim \ref{algoanalysisnum} is at most
$\exp_2\{O(g\log d/d + t \log(d^5/\gd)/\psi)\}$. If we replace
$d/2$ by $\psi$ in Step $1$ of the degree algorithm and take the
output of Step $2$ to be the final output, then the bound in the
claim becomes weaker:
$$
{2g\log d \choose \leq \frac{g}{\psi}}{2gd^3\log d \choose \leq
\frac{td}{(d-\psi)\psi}}=\exp_2\left\{O(g\log d/\psi + t
\log(d^5/\gd)/\psi)\right\}.
$$
(Each iteration of Step $1$ now reduces $G\setminus F$ by at least
$\psi$). Using this bound in the proof of Theorem \ref{thm-main}
instead of the stronger bound given by the three-step degree
algorithm would ultimately lead to a weaker bound on $\gb(\gl)$ in
(\ref{inq-gen.bd.on.lambda.1}). Step $1$ of the degree algorithm
may be though of as an ``initialization'' which reduces
$|G\setminus F|$ from $O(g)$ to $O(t)$ without adding much to the
``cost'' of the algorithm.
\end{remark}

\bigskip

\noindent {\em Proof of Lemma \ref{lem-reconstruction}: }Say that
$S$ is {\em small} if $|S| < g - \gamma t$ and {\em large}
otherwise. We can obtain all $A \in \cA$ for which $F \subseteq G$
and $S \supseteq [A]$ as follows.

If $S$ is small, we specify of $A$ by picking a subset of $S$. If
$S$ is large, we first specify $G$. Note that by
(\ref{inq-bounding.s.by.f}) and the definition of large we have in
this case that
$$
|G \setminus F| < 2t\psi/(d-\psi)+\gamma t~~~~\mbox{and}~~~~G
\setminus F \subseteq N(S) \setminus F,
$$
so we specify $G$ by picking a subset of $N(S) \setminus F$ of
size at most $2t\psi/(d-\psi)+\gamma t$ (this is our choice of $G
\setminus F)$. Then, noting that $[A]$ is determined by $G$, we
specify $A$ by picking a subset of $[A]$.

This procedure produces all possible $A$'s (and more). We now
bound the sum of the weights of the outputs.

If $S$ is small then the total weight of outputs is at most
\begin{equation} \label{ssmall}
(1+\gl)^{g - \gamma t}.
\end{equation}
We have
$$
|N(S)\setminus F| \leq d|S| \leq dg+\frac{2td\psi}{d-\psi} \leq
3dg
$$
so that if $S$ is large, the total number of possibilities for $|G
\setminus F|$ is at most
$$
{3dg \choose \leq \frac{2t\psi}{d-\psi}+\gamma t}
$$
and the total weight of outputs is at most
\begin{equation} \label{slarge}
{3dg \choose \leq \frac{2t\psi}{d-\psi} + \gamma t}(1+\gl)^{g-t}.
\end{equation}

Combining (\ref{ssmall}) and (\ref{slarge}), the lemma follows.
\qed

\medskip

\noindent {\it Acknowledgment}. This work originated while the
first author was a member of and while the second author was
visiting the Theory Group at Microsoft Research in Redmond,
Washington. The authors thank Microsoft Research for this support.

\end{document}